\documentclass[smallextended]{svjour3}
\usepackage{amsfonts}

\usepackage{latexsym}
\usepackage{amssymb}
\usepackage{amsmath}
\usepackage[mathscr]{eucal}
\usepackage{graphicx}
\usepackage{hyperref}
\usepackage{caption}
\usepackage{subcaption}

\renewcommand{\qed}{\hfill{\ \ \rule{2mm}{2mm}} \vspace{0.2in}}

\renewcommand{\thefigure}{\arabic{figure}}
\begin{document}

\title{Dirac-type Theorems for Inhomogenous Random Graphs}
\author{ \textbf{Ghurumuruhan Ganesan}
\thanks{E-Mail: \texttt{gganesan82@gmail.com} } \\
\ \\
IISER Bhopal, India}
\date{}
\maketitle


\begin{abstract}
In this paper, we study Dirac-type theorems for an inhomogenous  random graph~\(G\) whose edge probabilities are not necessarily all the same. We obtain sufficient conditions for the existence of Hamiltonian paths and perfect matchings, in terms of the \emph{sum} of edge probabilities. For edge probability assignments with two-sided bounds, we use P\'osa rotation and single vertex exclusion techniques to show that~\(G\) is Hamiltonian with high probability. For weaker one-sided bounds, we use bootstrapping techniques to obtain a perfect matching in~\(G,\) with high probability. We also highlight an application of our results in the context of channel assignment problem in wireless networks.

\vspace{0.1in} \noindent \textbf{Key words:} Dirac-type Theorems; Hamiltonian Paths; Perfect Matchings; Inhomogenous Random Graphs.

\vspace{0.1in} \noindent \textbf{AMS 2000 Subject Classification:} Primary: 05C62;
\end{abstract}

\bigskip

\renewcommand{\theequation}{\thesection.\arabic{equation}}
\setcounter{equation}{0}
\section{Introduction} \label{intro}
Hamiltonian cycles are important objects of study from both theoretical and application perspectives and for deterministic graphs, there are many known sufficient conditions for their existence. For example, Dirac~\cite{dirac} proved that if the minimum degree of a graph~\(H\) on~\(n\) vertices is at least~\(\frac{n}{2},\) then~\(H\) is Hamiltonian and for further results along this direction, we refer to Chapter~\(7\) of~\cite{west}.

In the case of random graphs, one of the earliest studies was initiated by~\cite{posa}, who used rotation techniques to obtain edge probability threshold for the occurrence of Hamiltonian paths with high probability, i.e. with probability converging to one as~\(n \rightarrow \infty.\) For a detailed survey of further techniques, we refer to the book~\cite{boll}. In~\cite{lee} a Dirac-type theorem is obtained for random subgraphs of a graph that is not necessarily complete. The main result there is that if the minimum degree of the ``parent" graph is at least~\(\frac{n}{2},\) then the Hamiltonicity threshold is again of the order of~\(\frac{\log{n}}{n}\) and they also demonstrated the tightness of this bound. We also remark that recently~\cite{cond} has obtained Dirac-type theorems random regular graphs.

In this paper, we consider inhomogenous random graphs with independently open edges, whose edge probabilities are not all the same.   We obtain sufficient ``two-sided" conditions on the edge probabilities, for the existence of a Hamiltonian path with high probability. We then show that under weaker one-sided conditions, there is a perfect matching with high probability.

The paper is organized as follows. In Section~\ref{sec_inhom} we state and prove our results regarding existence of Hamiltonian paths and perfect matchings in inhomogenous random graphs and in Section~\ref{sec_rand_env}, we briefly describe an application of our results in the context of channel assignment in wireless networks.

\renewcommand{\theequation}{\arabic{section}.\arabic{equation}}
\setcounter{equation}{0}
\section{Inhomogenous Random Graphs}\label{sec_inhom}
Let~\(K_n\) be the complete graph on~\(n\) vertices and let~\(\{Z(f)\}_{f \in K_n}\) be independent random variables with distribution
\begin{equation}\label{x_dist}
\mathbb{P}(Z(f) = 1) = p(f) = 1-\mathbb{P}(Z(f) = 0)
\end{equation}
where~\(0 \leq p(f) \leq 1\) and let~\(G\) be the random graph formed by the union of all edges~\(f\) satisfying~\(Z(f) = 1.\) A path in~\(G\) is a sequence of vertices~\(P = (u_1,\ldots,u_t)\) such that~\(u_i\) is adjacent to~\(u_{i+1}\) in~\(G,\) for each~\(1 \leq i \leq t-1.\)   The length of~\(P\) is the number of edges in~\(P\) which in this case is~\(t-1.\) If~\(P\) has the maximum possible length of~\(n-1,\) then we say that~\(P\) is a \emph{Hamiltonian} path.





Our first result obtains Hamiltonian paths in~\(G\) for the following class of ``well-behaved" probability assignments.  For constants~\(0 < \alpha < 1\) and~\(c_1,c_2 > 0\) and~\(0 < p = p(n) < 1,\) we say that the probability assignment~\(\{p(u,v)\}_{1 \leq u < v \leq n}\) in~(\ref{x_dist}) is~\((\alpha,c_1,c_2,p)-\)good if for any vertex~\(u\) and any set~\(S\) containing~\(r \geq \alpha n\) vertices,
\begin{equation}\label{good_cond}
c_1 rp \leq \sum_{v \in S} p(u,v) \leq c_2 rp.
\end{equation}
Here and throughout, constants do not depend on~\(n.\)  Letting~\(E_{HAM}\) be the event that~\(G\) contains a Hamiltonian path, we have the following result.
\begin{theorem}\label{lem_path}
Suppose~\(\{p(u,v)\}\) is~\((\alpha,c_1,c_2,p)-\)good for~\(p  = Cn^{-\frac{k}{k+1}}\) where\\\(C > 0\) is a constant and~\(k \geq 1\) is an integer constant. For every
\begin{equation}\label{th_max}
C \leq \frac{1}{10c_2} \text{ and }\alpha \leq \min\left(\frac{7}{8},\left(\frac{Cc_1}{8}\right)^{k+1}\right),
\end{equation}
there exists a constant~\(\theta > 0\) such that
\begin{equation}\label{ham_bound}
\mathbb{P}(E_{HAM}) \geq 1-  e^{-\theta np}.
\end{equation}
\end{theorem}
In other words, the resulting random graph~\(G\) contains a Hamiltonian path with high probability, i.e., with probability converging to one as~\(n \rightarrow \infty.\)

For example, suppose we consider the particular case when~\(G\) is homogenous with common edge probability~\(p_{hom} \geq \frac{1}{n^{\beta}}\) for some~\(0 < \beta < 1.\) If~\(k\) is a large enough integer so that~\(0 < \beta < \frac{k}{k+1} < 1,\) then by direct coupling, it suffices to demonstrate Hamiltonicity for the homogenous random subgraph~\(G_{low}\) with edge probability~\(p = Cn^{-\frac{k}{k+1}}\) where~\(C\) is chosen sufficiently small so that the first condition in~(\ref{th_max}) holds. In this case condition~(\ref{good_cond}) holds trivially with~\(c_1 = c_2=1\)  and all constant~\(\alpha >0.\) Therefore the second condition in~(\ref{th_max}) is also satisfied and so from~(\ref{ham_bound}), we then get that~\(G_{low}\) (and therefore~\(G\)) contains a Hamiltonian path with high probability.

Our proof of Theorem~\ref{lem_path} involves a combination of the P\'osa rotation technique~\cite{posa} and the ``single vertex exclusion" method used in~\cite{ganesan} for obtaining Hamiltonicity in inhomogenous random graphs when~\(k = 1.\) Throughout, we use the following standard deviation estimate. Let~\(Z_i, 1 \leq i \leq t\) be independent Bernoulli random variables satisfying~\(\mathbb{P}(Z_i = 1) = p_i = 1-\mathbb{P}(Z_i = 0).\) If~\(W_t = \sum_{i=1}^{t} Z_i\) and~\(\mu_t = \mathbb{E}W_t,\) then for any~\(0 < \eta < \frac{1}{2}\) we have that
\begin{equation}\label{conc_est_f}
\mathbb{P}\left(\left|W_t-\mu_t\right| \geq \eta \mu_t\right) \leq 2\exp\left(-\frac{\eta^2}{4}\mu_t\right).
\end{equation}
For a proof of~(\ref{conc_est_f}), we refer to Corollary~\(A.1.14,\) pp.~\(312,\) Alon and Spencer (2008).\\\\
\emph{Proof of Theorem~\ref{lem_path}}: We obtain the  Hamiltonian path of~\(G\) in three steps with the first two steps are preliminary calculations. In the first step, we define an event~\(E_{good}\) regarding the neighbourhood size of arbitrary sets and obtain probability estimates for~\(E_{good}.\) Next, in the second step, we define the concept of pivots and use the occurrence of~\(E_{good}\) to estimate the number of pivots in the maximum length path of~\(G\) obtained after excluding a single vertex. Finally, in the third step, we combine the two estimates to show that~\(\pi\) is Hamiltonian with high probability.

\underline{\emph{Step 1}}: Let~\(S \subset \{1,2,\ldots,n\} \) be any set of size~\(s  := \#S\) and let~\({\cal N}_{out}(S)\) be the set of all vertices of the complement set~\(S^c,\) adjacent to at least one vertex of~\(S.\) Letting~\(N_{out}(S) = \#{\cal N}_{out}(S)\) denote the size of~\({\cal N}_{out}(S),\) we begin by showing that for any set~\(S\) of size~\(1 \leq s \leq \frac{1}{10c_2p},\) we have
\begin{equation}\label{out_cond}
\frac{3c_1nps}{4}  \leq \mathbb{E} N_{out}(S) \leq c_2nps.
\end{equation}
Indeed, by definition~\(s \leq \frac{1}{10c_2p}=o(n)\) and so~\(n-s \geq \frac{7n}{8}\) for all~\(n\) large. Therefore using~(\ref{good_cond}), we get that the expected number of vertices adjacent to~\(u\) in~\(G\) is at least~\(c_1np\) and at most~\(c_2np.\) This proves~(\ref{out_cond}) for the case~\(s=1\) and also the upper bound in~(\ref{out_cond}) for general~\(s.\)

We prove the lower bound in~(\ref{out_cond}) by induction on~\(s.\) Pick a set~\(S\) of size~\(2 \leq s \leq \frac{1}{10c_2p}\) and for~\(u \in S,\) define~\(S_u := S \setminus \{u\}.\) The set~\(S_u\) has size~\(s-1\) and so by the induction assumption
\begin{equation}\label{gisele_tits}
\frac{3c_1np(s-1)}{4} \leq \mathbb{E} N_{out}(S_u) \leq c_2np(s-1).
\end{equation}
Using the concentration estimate~(\ref{conc_est_f}) with~\(\eta >0\) small, we then get that
\begin{equation}\label{zilpa}
\mathbb{P}\left(N_{out}(S_u) \leq (1+\eta)c_2np(s-1)\right) \geq 1-e^{-C_0np(s-1)} \geq 1-e^{-C_0np}
\end{equation}
for some constant~\(C_0 > 0,\) not depending on~\(s.\) If the event defined in the left hand side of~(\ref{zilpa}) occurs, then there are at least~\(n-(1+\eta)c_2np(s-1)\) vertices in~\(S_u^c\) that are not adjacent to any vertex of~\(S_u\) and at least~\(n-1-(1+\eta)c_2np(s-1)\) among these are adjacent to~\(u\) in~\(K_n.\)

Recalling that~\(s \leq \frac{1}{10c_2p},\) we then choose~\(\eta > 0\) small so that~\[n-1-(1+\eta)c_2np(s-1) \geq n-1-\frac{(1+\eta)}{10} \geq \frac{7n}{8}\] for all~\(n \geq N\) large not depending on~\(s.\) With this choice of~\(\eta\) we get from~(\ref{good_cond}) that the expected number of vertices of~\(S^c\) adjacent \emph{only} to~\(u\) and no other vertex of~\(S,\) is at least~\(\frac{7c_1np}{8}(1-e^{-C_0np}).\) Thus
\[\mathbb{E}N_{out}(S) \geq  \frac{3c_1np(s-1)}{4} + \frac{7c_1np}{8}(1-e^{-C_0np}) \geq \frac{3c_1 nps}{4}\] and this proves the induction step for the lower bound in~(\ref{out_cond}).

Defining~\(E_{good}(S) := \left\{\frac{c_1nps}{2} \leq N_{out}(S) \leq 4c_2nps\right\}\) and using the deviation estimate~(\ref{conc_est_f}) with~\(\eta = \frac{1}{4},\) we obtain from the bounds in~(\ref{out_cond}) that
\[\mathbb{P}(E_{good}(S)) \geq 1-\exp\left(-2C_1 nps\right)\] for some constant~\(C_1 > 0.\) There are at most~\({n \choose s} \leq n^{s}\) choices for~\(S\) and so setting~\(E_{good} := \bigcap_{S} E_{good}(S),\) where the intersection is over all permissible~\(S,\) we get from the union bound that
\[\mathbb{P}(E_{good}) \geq 1- \sum_{s \geq 1} n^{s}\exp\left(-2C_1 nps \right).\] Since~\(p = Cn^{-\frac{k}{k+1}}\) is much larger than~\(\frac{\log{n}}{n},\) we have that
\[n^{s}\exp\left(-2C_1 nps \right) = \exp\left(s\left(\log{n} - 2C_1 np \right)\right) \leq \exp\left(-C_1 nps \right),\]
and so
\begin{equation}\label{e_good_est}
\mathbb{P}(E_{good}) \geq 1-\sum_{s  \geq 1} e^{-C_1 nps} \geq 1-2e^{-C_1np}
\end{equation}
for all~\(n\) large.

\begin{figure}[tbp]
\centering
\includegraphics[width=1.5in, trim= 100 340 100 100, clip=true]{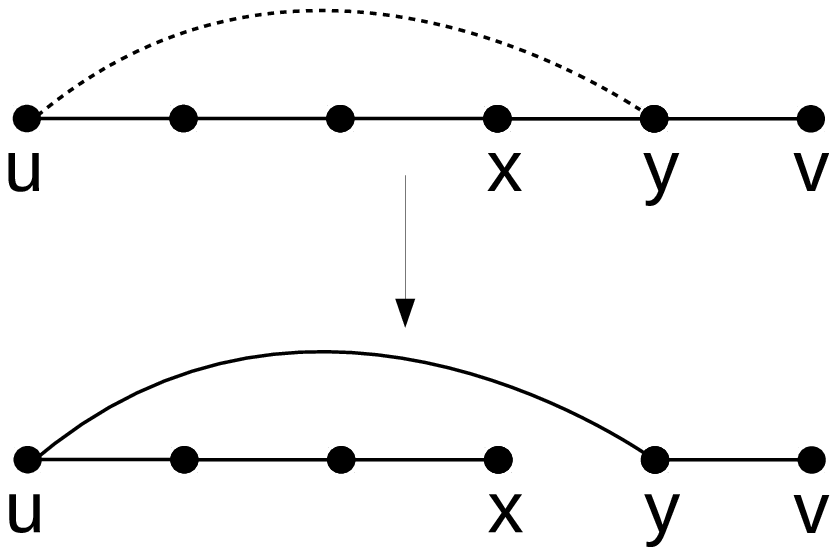}
\caption{The path~\(\pi\) with endvertices~\(u\) and~\(v\) is shown at the top, together with a neighbour~\(y\) of~\(u.\) We perform a P\'osa rotation to obtain a new path with the pivot~\(x\) as the endvertex, as illustrated in the bottom figure.}
\label{fig_ham}
\end{figure}

The event~\(E_{good}\) as defined above, allows us to count the number of pivots in the maximum length paths of~\(G\) as demonstrated in Step~\(2\) below.\\\\
\underline{\emph{Step 2}}: Let~\(\pi = (\pi(1),\pi(2),\ldots,\pi(t))\) be the longest path in~\(G\) with endvertices~\(\pi(1)\) and~\(\pi(t).\) If~\({\cal S}_1\) is the set of neighbours of~\(\pi(1)\) in~\(G,\) then all vertices in~\({\cal S}_1\) must also be present in~\(\pi\) and so~\({\cal S}_1 = \{\pi(a_1),\ldots,\pi(a_r)\}\) for some integers~\(2 = a_1  < a_2 < \ldots < a_r.\) For each~\(a_j, j \geq 1\) we can do a P\'osa rotation as shown in Figure~\ref{fig_ham} and obtain a maximum length path with endvertex~\(\pi(a_{j}-1).\) We refer to~\(\pi(a_j-1), 2 \leq j \leq r\) as \emph{first generation pivots} or pivots associated to~\(\pi(1).\)

Consider now the pivot~\(\pi(a_j-1)\) and let~\(\gamma_j = (\gamma_j(1),\ldots,\gamma_j(t))\) be a maximum length path in~\(G\) with endvertex~\(\gamma_j(1) = \pi(a_j-1),\) obtained via the P\'osa rotation as described in the previous paragraph. Again all neighbours of~\(\pi(a_j-1)\) in the graph~\(G\) must be present in~\(\gamma_j\) as well and if~\({\cal S}_j\) is the set of neighbours of~\(\pi(a_j-1)\) in~\(G,\) then we obtain~\(\#{\cal S}_j\) pivots associated to~\(\pi(a_j-1),\) which we call as \emph{second generation pivots}. Continuing this way we define~\(i^{th}-\)generation pivots for arbitrary~\(i \geq 1.\)  Summarizing, the longest path~\(\pi\) in~\(G\) satisfies the following properties:\\
\((p1)\) All neighbours of the two endvertices of~\(\pi\) are present in~\(\pi\) itself.\\
\((p2)\) If~\(v\) is a~\(i^{th}-\)generation pivot of~\(\pi\) and~\(\gamma_v\) is a maximum length path obtained from successive P\'osa rotations of~\(\pi\) and containing~\(v\) as an endvertex, then all neighbours of~\(v\) belong to~\(\gamma_v\) (and therefore to~\(\pi\)) as well.

We now assume that the event~\(E_{good}\) described in Step~\(1\) occurs and estimate the number of pivots in the maximum length path obtained by excluding a single vertex. Specifically, we perform an iterative pivot counting procedure consisting of~\(k\) steps at the end of which we demonstrate that the overall number of pivots grows  linearly in~\(n.\) Details follow.

For~\(1 \leq j \leq n\) we let~\(\pi_j\) be the maximum length path in the graph~\(G_j\) obtained by removing the vertex~\(j\) from~\(G\) and show by induction that  for each~\(1 \leq l \leq k,\) there are at least~\(\left(\frac{c_1np}{8}\right)^{l}\) and at most~\((8c_2np)^{l}\) pivots belonging to the~\(l^{th}\) generation (here we recall that the integer~\(k\) is defined via~\(p = Cn^{-\frac{k}{k+1}}\)); i.e., our goal is to obtain the bounds
\begin{equation}\label{pl_est}
\left(\frac{c_1 np}{8}\right)^{l} \leq \#{\cal P}_l \leq (8c_2np)^{l}
\end{equation}
for each~\(1 \leq l \leq k.\)

We begin with the base case~\(l=1.\) Indeed if~\(\pi_j = (\pi_j(1),\ldots,\pi_j(t))\) then since~\(E_{good}\) occurs, the vertex~\(\pi_j(1)\) contains at least~\(\frac{c_1np}{2}\) neighbours and at most~\(4c_2np\) neighbours in~\(G.\) Therefore there are at least~\(\frac{c_1np}{2}-1 \geq \frac{c_1np}{8}\) vertices adjacent to~\(\pi_j(1)\) in~\(G_j\) and so there are at least~\(\frac{c_1np}{8}\) pivots associated with~\(\pi_j(1).\) Each vertex adjacent to~\(\pi_j(1)\) in~\(G_j\) gives rise to at most two pivots and so there are at least~\(\frac{c_1np}{8}\) and at most~\(2(4c_2np) = 8c_2np,\) first generation pivots. This completes the proof of~(\ref{pl_est}) for the case~\(l=1.\)

To argue for general~\(l,\) we first see that the~\(E_{good}\) probability estimate is valid as long as~\((8c_2np)^{l}  \leq \frac{1}{10c_2p}\) (see step~\(1\)). Since~\(p = \frac{C}{n^{k/k+1}},\) we have that~\[p(8c_2np)^{k} = C^{k+1}(8c_2)^{k} \leq \frac{1}{10c_2}\] if~\(C > 0\) satisfies the first condition in~(\ref{th_max}). Fixing such a~\(C\) henceforth, we now proceed to the induction step.  For~\(l \geq 1,\) let~\({\cal P}_l\) be the set of all~\(l^{th}-\)generation pivots of~\(\pi_j\) and let~\(\bigcup_{v \in {\cal P}_l} {\cal N}_j(v)\) be the set of all vertices adjacent to at least one vertex of~\({\cal P}_l\) in~\(G_j,\) where~\({\cal N}_j(v)\) denotes the set of all neighbours of~\(v\) in~\(G_j.\)   By induction assumption
\begin{equation}\label{piv_est}
\left(\frac{c_1 np}{8}\right)^{l} \leq \#{\cal P}_l \leq (8c_2np)^{l}
\end{equation} and due to the occurrence of the event~\(E_{good},\) we have that
\begin{equation}\label{tot_nei_est}
\frac{c_1np}{2} \cdot \left(\frac{c_1np}{8}\right)^{l}-1 \leq \#\left(\bigcup_{v \in {\cal P}_l} {\cal N}_j(v)\right) \leq (4c_2np) (8c_2np)^{l}.
\end{equation}
Thus the number of~\((l+1)^{th}-\)generation pivots is at most~\[2\cdot(4c_2np) \cdot (8c_2np)^{l} = (8c_2np)^{l+1}.\]

For a lower bound on the number of~\((l+1)^{th}-\)generation pivots, we pick a~\(l^{th}-\)generation pivot~\(v \in {\cal P}_l\) and perform~\(l\) P\'osa rotations to obtain a path~\(\gamma_{v,j}\) containing~\(v\) as an endvertex. In doing so, we remove exactly~\(l\) edges from~\(\pi_j\) and add~\(l\) other edges, to obtain~\(\gamma_{v,j}.\)  We let~\({\cal R}_{j}(v)\) be the union of the set of all  endvertices of the~\(l\) removed edges,~\(l\) newly added edges and the endvertices of~\(\gamma_{v,j}\) so that there are at most~\(4l+2\) vertices in~\({\cal R}_j(v).\)

We now obtain a lower bound on the number of pivots that arise out of \emph{some}~\(l^{th}\) generation pivot in~\( \Lambda_j := \bigcup_{v \in {\cal P}_l} \left({\cal N}_j(v) \setminus {\cal R}_j(v)\right).\) Let~\(w \in \Lambda_j\) be adjacent to~\(v \in {\cal P}_l\) in~\(G_j.\) From property~\((p2)\) described above, all neighbours of~\(v\) in~\(G_j\) belong to the path~\(\pi_j\) and so~\(w = \pi_j(x)\) for some~\(x.\) Moreover, because~\(w \notin {\cal R}_j(v),\) both the edges~\[(\pi_j(x),\pi_j(x+1)) \text{ and }(\pi_j(x-1),\pi_j(x))\] must belong to the path~\(\gamma_{v,j}\) with~\(v\) as an endvertex, obtained by performing~\(l\) P\'osa rotations on~\(\pi_j.\) Thus at least one of the vertices in \[ \{\pi_j(x-1),\pi_j(x+1)\},\] say for example~\(\pi_j(x-1),\) must necessarily be a~\((l+1)^{th}-\)generation pivot and we say that~\(\pi_j(x-1)\) is an~\((l+1)^{th}-\)generation pivot \emph{created} by an~\(l^{th}-\)generation pivot.

Conversely, any~\((l+1)^{th}-\)generation pivot created from an~\(l^{th}-\)generation pivot is of the form~\(\pi_j(y)\) and is created from an~\(l^{th}-\)generation pivot that is necessarily adjacent to either~\(\pi_j(y-1) \in \Lambda_j\) or~\(\pi_j(y+1) \in \Lambda_j.\) We have a couple of remarks here. It may happen that more than one~\(l^{th}-\)generation pivot itself could be adjacent to either~\(\pi_j(y-1)\) or~\(\pi_j(y+1)\) or some~\((l+1)^{th}-\)generation pivot is already a~\(y^{th}-\)generation pivot for some~\( y \leq l.\) In any case, our count of the number of pivots in the~\((l+1)^{th}\) generation above depends only on the \emph{neighbours} of the~\(l^{th}-\)generation pivots, i.e., vertices in~\(\Lambda_j.\)  Therefore from the above argument, we get
\begin{align}
\#{\cal P}_{l+1} &\geq \frac{\#{\Lambda_j}}{2} \nonumber\\
&= \frac{1}{2}\#\left(\bigcup_{v \in {\cal P}_l} {\cal N}_j(v) \setminus {\cal R}_j(v)\right) \nonumber\\
&\geq \frac{1}{2} \#\left( \left(\bigcup_{v \in {\cal P}_l} {\cal N}_j(v)\right) \setminus \left(\bigcup_{v \in {\cal P}_l} {\cal R}_j(v)\right)\right) \nonumber\\
&\geq \frac{1}{2}\#\left(\bigcup_{v \in {\cal P}_l} {\cal N}_j(v)\right) - \sum_{v \in {\cal P}_l} \#{\cal R}_j(v)\nonumber\\
&\geq \frac{1}{2}\#\left(\bigcup_{v \in {\cal P}_l} {\cal N}_j(v)\right) -(4l+2)\#{\cal P}_l, \label{mgr}
\end{align}
since the number of vertices in any~\({\cal R}_j(v)\) is at most~\(4l+2.\) Plugging the bounds from~(\ref{piv_est}) and~(\ref{tot_nei_est}) into~(\ref{mgr}), we get that
\[\#{\cal P}_{l+1} \geq \frac{1}{2}\left(\frac{c_1np}{2} \cdot \left(\frac{c_1np}{8}\right)^{l}-1\right) -  (4l+2)(8c_2np)^{l} \geq \left(\frac{c_1np}{8}\right)^{l+1}\] for all~\(n\) large. This obtains~(\ref{pl_est}) for the general case.

From~(\ref{pl_est}), we see that the number of~\((k+1)^{th}-\)generation pivots is at least~\[\left(\frac{c_1np}{8}\right)^{k+1} = \left(\frac{c_1}{8}\right)^{k+1} \cdot n \cdot p(np)^{k} =: Dn,\] where~\(D  = \left(\frac{Cc_1}{8}\right)^{k+1}>0\) is a constant. Summarizing, we have shown that if the event~\(E_{good}\) occurs then there are at least~\(Dn\) pivots in~\(\pi_j\) for any~\(1 \leq j \leq n\) and this completes the second step of the proof.



\underline{\emph{Step 3}}: We now combine the estimates obtained in the previous two steps to obtain the desired Hamiltonian path. We begin by noting that if the vertex~\(j\) does not belong to the maximum length path~\(\pi,\) then the longest path~\(\pi_j\) in the random graph~\(G_j\) is also~\(\pi\) and more importantly,~\(j\) is \emph{not} adjacent to any pivot of~\(\pi_j = \pi;\) otherwise, we would have a path containing~\(j\) as an endvertex and with length strictly larger than~\(\pi.\) Therefore letting~\(Q_j(\pi_j)\) be the event that~\(j\) is not adjacent to any pivot of~\(\pi_j,\) we get that \begin{equation}\label{thnoof}
\mathbb{P}\left(\{j \notin \pi\} \cap E_{good}\right)  \leq \sum_{\Gamma} \mathbb{P}\left(\{\pi_j = \Gamma\} \cap Q_j(\Gamma)\right)
\end{equation}
where the summation is over all deterministic paths~\(\Gamma\) not containing~\(j,\) satisfying~\((p1)-(p2)\) and containing at least~\(Dn\) pivots. The event~\(Q_j(\Gamma)\) is independent of the event~\(\{\pi_j = \Gamma\}\) and moreover if the probability assignment is~\((\alpha,c_1,c_2,p)-\)good for~\(\alpha \leq D,\) then we get that~\(\mathbb{P}(Q_j(\Gamma)) \leq e^{-c_1Dnp}.\) Consequently,
\begin{eqnarray}
\mathbb{P}\left(\{\pi_j = \Gamma\} \cap Q_j(\Gamma)\right) &=& \mathbb{P}(\pi_j = \Gamma) \mathbb{P}(Q_j(\Gamma)) \nonumber\\
&\leq& e^{-c_1Dnp}\mathbb{P}(\pi_j = \Gamma) \nonumber
\end{eqnarray}
Summing over~\(\Gamma\) we therefore get from~(\ref{thnoof}) that~\(\mathbb{P}\left(\{j \notin \pi\} \cap E_{good}\right) \leq   e^{-c_1Dnp}\)
and so using~(\ref{e_good_est}), we get that
\[\mathbb{P}(j \notin \pi) \leq e^{-c_1Dnp} + \mathbb{P}(E^c_{good}) \leq e^{-c_1Dnp} + 2e^{-C_1np} \leq e^{-C_2np}\] for some constant~\(C_2 > 0.\)
Therefore by the union bound, we get that~\(\pi\) is Hamiltonian with probability at least~\(1-ne^{-C_2 np}\) and this completes the proof of Theorem~\ref{lem_path}.~\(\qed\)



In Theorem~\ref{lem_path} we have obtained sufficient conditions to be satisfied by the probability assignments, for the existence of Hamiltonian paths.
Specifically,~(\ref{good_cond}) has ``two-sided" conditions on both the upper and lower bounds of the edge probability sums. A natural question therefore is whether under weaker conditions, we can say anything about the existence of Hamiltonian paths or perhaps at least perfect matchings. For context, we recall that a set of vertex disjoint edges in the random graph~\(G\) is called a \emph{matching} and a matching~\({\cal W}\) of~\(G\) is said to be \emph{perfect} if every vertex except at most one, is an endvertex of some edge in~\({\cal W}.\) By definition, if~\(G\) has a Hamiltonian path~\(\pi,\) then picking every alternate edge of~\(\pi\) gives us a perfect matching. Thus the existence of a Hamiltonian path implies the existence of a perfect matching but the converse, however, is not true.

We now consider probability assignments with ``one-sided" conditions involving only lower bounds and obtain perfect matchings in the resulting random graph~\(G.\) For constants~\(0 < \beta,\gamma < 1\) and~\(0 < p = p(n) < 1,\) we say that~\(\{p(u,v)\}\) is a~\((\beta,d_1,d_2,p)-\)nice probability assignment if the following holds for any two vertices~\(u,v:\) For any two disjoint sets~\(S_1,S_2\) not containing~\(u\) or~\(v\) and having cardinality~\(r \geq \beta n\) each, we have
\begin{equation}\label{goo_cond2}
\sum_{w \in S_1} p(u,w) \geq d_1 rp\;\;\text{ and } \sum_{w_1 \in S_1,w_2 \in S_2}p(u,w_1)p(w_2,v) \geq d_2 rp^2.
\end{equation}
Letting~\(E_{PER}\) be the event that~\(G\) contains a perfect matching, we have the following result. As before, constants do not depend on~\(n.\)
\begin{theorem}\label{thm_match} Suppose~\(\{p(u,v)\}\) is~\((\beta,d_1,d_2,p)-\)nice for some~\(p  \geq \frac{\log{n}}{\sqrt{n}}\) and\\\(\beta \leq \frac{1}{4}.\) There is a constant~\(D  > 0\) such that
\begin{equation}\label{per_bound}
\mathbb{P}(E_{PER}) \geq 1-  e^{-D(\log{n})^2}.
\end{equation}
\end{theorem}
In other words, the one-sided conditions in~(\ref{goo_cond2}) ensure the existence of perfect matchings with high probability. In the next Section, we briefly describe an application of our result in the context of channel assignment in wireless networks.


\emph{Proof of Theorem~\ref{thm_match}}: For simplicity we assume throughout that~\(n = 4z\) is even and let~\(G_{bip}\) be the bipartite subgraph of~\(G\) with left vertex set~\(X = \{1,2,\ldots,2z\}\) and right vertex set~\(Y = \{2z+1,\ldots,n\}.\) We obtain a perfect matching in~\(G_{bip}\) using a bootstrapping technique as follows. In the first step, we let~\({\cal M}\) be the maximum matching in~\(G_{bip}\) (picked according to a deterministic rule) and estimate the probability~\(p_{0}\) that~\(\#{\cal M} \geq z.\) In the next step, we then use the estimate for~\(p_{0}\) as a bootstrap to bound the probability of a perfect matching.  We provide the details below.

For~\({\cal L} \subseteq X\) and~\({\cal R} \subseteq Y\) let~\(E\left({\cal L}, {\cal R}\right)\) be the event that no vertex of~\({\cal L}\) is adjacent to any vertex of~\({\cal R},\) in~\(G_{bip}.\) If~\(\#{\cal L} \geq z\) and~\(\#{\cal R} \geq  z,\) then any vertex~\(v \in {\cal L}\) is adjacent to at least~\( z = \frac{n}{4}\) vertices of~\({\cal R}\) in~\(K_n\) and so using the first condition in~(\ref{goo_cond2}), we get that the expected number of vertices adjacent to~\(v\) in~\(G_{bip}\) is at least~\(d_1 zp.\) Therefore the \emph{total} expected number of edges of~\(G_{bip}\) containing one endvertex in~\({\cal L}\) and the other endvertex in~\({\cal R}\) is at least~\(d_1 z^2p = \frac{d_1 n^2p}{16}\) and consequently,
\begin{equation}\label{eq_one}
\mathbb{P}\left(E\left({\cal L}, {\cal R}\right)\right) \leq \exp\left(-\frac{d_1  n^2p}{16}\right).
\end{equation}

Now if the maximum matching~\({\cal M}\) of~\(G_{bip}\) has size~\(\#{\cal M} <  z,\) then removing the edges and vertices of~\({\cal M}\) from~\(G_{bip}\) we are left with an edgeless bipartite graph containing at least~\(z\) left vertices and at least~\(n-2z-z \geq z\) right vertices. Thus~\(\mathbb{P}\left(\#{\cal M} < z\right) \leq \mathbb{P}\left(\bigcup E\left({\cal L}, {\cal R}\right)\right)\) where the union is over all sets~\({\cal L} \subseteq X\) and~\({\cal R} \subseteq Y\) satisfying~\(\#{\cal L} \geq  z\) and~\(\#{\cal R}\geq  z.\) Letting~\(E_{low} := \{\#{\cal M} \geq  z\},\) we then get from the union bound and~(\ref{eq_one}) that
\begin{equation}\label{e_low_est}
\mathbb{P}\left(E_{low}^c\right) \leq \sum_{{\cal L}, {\cal R}} \mathbb{P}\left(E\left({\cal L},{\cal R}\right)\right) \leq 4^{n}\cdot \exp\left(-\frac{d_1n^2p}{16}\right)
\end{equation}
since there are at most~\(2^{n}\) choices each, for~\({\cal L}\) and for~\({\cal R}.\)

\begin{figure}[tbp]
\centering
\includegraphics[width=1.5in, trim= 80 370 150 50, clip=true]{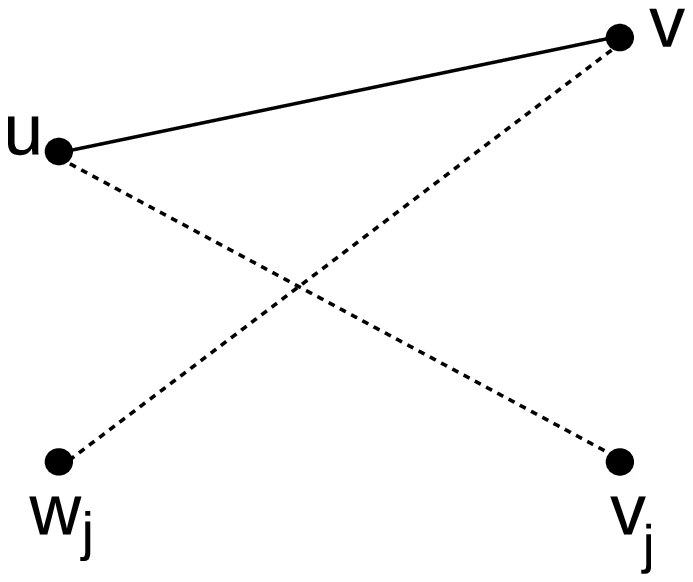}
\caption{Replacing the edge~\((u,v)\) with the edges~\((u,v_j)\) and~\((w_j,v)\) to get a bigger matching.}
\label{fig_match}
\end{figure}

Suppose the event~\(E_{low}\) occurs and for~\(1 \leq  j \leq z\) let~\(E_j\) be the event that~\({\cal M}\) does not have any edge with an endvertex in~\(\{w_j,v_j\}\) where~\( w_j \in X\) and~\(v_j \in Y.\) If~\(E_j \cap E_{low}\) occurs, then in the random graph~\(G^{(j)}\) obtained by removing the vertices~\(w_j\) and~\(v_j\) from~\(G_{bip},\) the maximum matching~\({\cal M}(G^{(j)})\) in~\(G^{(j)}\) is still~\({\cal M}.\)
Consequently
\begin{eqnarray}
\mathbb{P}(E_j \cap E_{low}) &=& \mathbb{P}\left(E_j \cap E_{low} \cap \{{\cal M}(G^{(j)}) = {\cal M}\}\right) \nonumber\\
&=& \sum_{{\cal E}} \mathbb{P}\left(E_j \cap \{{\cal M}(G^{(j)}) = {\cal E}\} \cap \{{\cal M} = {\cal E}\}\right)\label{ej_est_one}
\end{eqnarray}
where the summation in~(\ref{ej_est_one}) is over all sets of edges~\({\cal E}\) of size at least~\(z\) and satisfying the property that no edge of~\({\cal E}\) has an endvertex in the set~\(\{w_j,v_j\}.\) Moreover,~\[E_j \cap \{{\cal M}(G^{(j)}) = {\cal E}\} \cap \{{\cal M} = {\cal E}\} \subset \{{\cal M}(G^{(j)}) = {\cal E}\}  \cap V\left({\cal E}\right),\] where~\(V\left({\cal E}\right)\) is the event that there is no edge~\(e = (u,v) \in {\cal E}\) such that~\(v_j\) is adjacent to~\(u \in X\) and~\(w_j\) is adjacent to~\(v \in Y,\) in the graph~\(G_{bip}.\) This is because if there existed such an edge~\((u,v),\) then we could remove~\((u,v)\) from~\({\cal M}\) and add the edges~\((u,v_j)\) and~\((w_j, v)\) to get a matching of~\(G_{bip}\) whose size is strictly larger than that of~\({\cal E},\) a contradiction.  This is illustrated in Figure~\ref{fig_match}.

From~(\ref{ej_est_one}) and the discussion in the above paragraph, we therefore get that
\begin{eqnarray}
\mathbb{P}(E_j \cap E_{low}) &\leq& \sum_{{\cal E}} \mathbb{P}\left(\{{\cal M}(G^{(j)}) = {\cal E}\} \bigcap V\left({\cal E}\right)\right) \nonumber\\
&=& \sum_{{\cal E}} \mathbb{P}\left({\cal M}(G^{(j)}) = {\cal E}\right) \mathbb{P}\left(V\left({\cal E}\right)\right) \label{ej_est_two}
\end{eqnarray}
because for any~\({\cal E},\) the event~\(V\left({\cal E}\right)\) depends only on the state of edges containing an endvertex in~\(\{w_j,v_j\}\) and is therefore independent of~\(\{{\cal M}(G^{(j)}) = {\cal E}\},\) by definition. To estimate~\(\mathbb{P}(V({\cal E})),\) we let~\({\cal E} = \{(x_i,y_i)\},1 \leq i \leq t, t \geq z = \frac{n}{4},\) be the edges in~\({\cal E}\) and obtain that \[\mathbb{P}(V({\cal E})) \leq \prod_{i=1}^{t}\left(1-p(w_j,y_i)p(x_i,v_j)\right) \leq \exp\left(-\sum_{i=1}^{t}p(w_j,y_i)p(x_i,v_j)\right).\]  From the second condition in~(\ref{goo_cond2}) we get that~\[\mathbb{P}\left(V\left({\cal E}\right)\right) \leq e^{-d_2tp^2} \leq \exp\left(-\frac{d_2p^2n}{4}\right).\] Plugging this into~(\ref{ej_est_two}) we get
\begin{eqnarray}
\mathbb{P}(E_j \cap E_{low}) &\leq& \exp\left(-\frac{d_2p^2n}{4}\right)\sum_{{\cal E}} \mathbb{P}\left({\cal M}(G^{(j)}) = {\cal E}\right) \nonumber\\
&\leq& \exp\left(-\frac{d_2p^2n}{4}\right) \nonumber
\end{eqnarray}
and so from~(\ref{e_low_est}) we therefore get that
\[\mathbb{P}(E_j) \leq 4^{n}\cdot \exp\left(-\frac{d_1 n^2p}{16}\right)+ \exp\left(-\frac{d_2p^2n}{4}\right) =:Q.\]


Summarizing, we get that at least one of vertices~\(w_j\) and~\(v_j\) belong to the maximum matching of~\(G_{bip}\) with probability at least~\(1-Q.\) There are at most~\(n^2\) edges with one endvertex in~\(X\) and other endvertex in~\(Y\) and so applying the union bound, we see that~\(G_{bip}\) contains a perfect matching with probability at least~\(1-n^2Q.\) Using the fact that~\(p \geq \frac{\log{n}}{\sqrt{n}},\) we then get that~\(\mathbb{P}(E_{PER}) \geq 1-e^{-D_1(\log{n})^2}\) for some constant~\(D_1 > 0\) and this completes the proof of Theorem~\ref{thm_match}.~\(\qed\)

\renewcommand{\theequation}{\arabic{section}.\arabic{equation}}
\setcounter{equation}{0}
\section{Channel Assignment Problem}\label{sec_rand_env}
In this section, we briefly describe an application of our results in the context of channel assignment problem in wireless networks. Let~\(K_{n,n}\) be the complete bipartite graph with vertex sets~\(X = Y=\{1,2,\ldots,n\}\) and let~\(\{Z(u,v)\}_{1 \leq u,v \leq n} \) be positive independent random variables with~\(F_{u,v}\) denoting the distribution of~\(Z_{u,v}.\)

In the context of wireless networks, the set~\(X\) denotes the set of users and~\(Y\) denotes the set of channels to be assigned to the users with the constraint that no two users are assigned the same channel. The random variable~\(Z(u,v)\) denotes the fading gain~\cite{goldsmith} of the~\(u^{th}\) user on the~\(v^{th}\) channel and in many practical scenarios, the fading gains are independent but not necessarily identically distributed. It is of interest to assign the ``best possible" channel to each user and one straightforward way to implement this would be to set a predetermined threshold~\(\lambda\) and assign each user a distinct channel whose fading gain is at least~\(\lambda.\) The natural question is whether such an assignment does in fact exist and we use perfect matchings described in the previous section to demonstrate an answer.

Indeed, let~\(G_{bip}\) be the random bipartite graph obtained by retaining all edges~\((u,v)\) satisfying~\(Z(u,v) > \lambda\) and we set~\(p(u,v) := \mathbb{P}(Z_{u,v} > \lambda).\) If the condition~(\ref{goo_cond2}) in Theorem~\ref{thm_match} holds, then we are guaranteed a perfect matching in~\(G_{bip}\) with high probability, i.e., with probability converging to one as~\(n \rightarrow \infty.\) Any perfect matching in~\(G_{bip}\) provides a valid channel assignment as described before and so the condition~(\ref{goo_cond2}) could be interpreted as a sufficient condition for assigning each user  a distinct channel, with fading gain at least~\(\lambda.\)

Finally, we remark that the iterative analysis described in our paper could also be potentially extended to preferential attachment models~\cite{piva} and we plan to analyze applicability to these models in the future.

\subsection*{\em Data Availability Statement}
Data sharing not applicable to this article as no datasets were generated or analysed during the current study.

\subsection*{\em Acknowledgement}
I thank Professors Rahul Roy, Federico Camia, Alberto Gandolfi, C. R. Subramanian and the referee for crucial comments that led to an improvement of the paper. I also thank IMSc and IISER Bhopal for my fellowships.

\subsection*{\em Conflict of Interest and Funding Statement}
I certify that there is no actual or potential conflict of interest in relation to this article. No funds, grants or other support was received for the preparation of this manuscript.

\bibliographystyle{plain}

\end{document}